\documentclass[11pt]{article}
\usepackage{amsfonts}
\usepackage{amsmath}
      \def\FF{\mathbb F}
      
  \def\QQ{\mathbb Q}    
\def\ZZ{\mathbb Z}
\newcommand\scrA{{\mathcal{A}}}

\newcommand\End{\mathop{\rm End}\nolimits}
\newcommand\Aut{\mathop{\rm Aut}\nolimits}

\newcommand\GL{{\mathrm{GL}}}

\newcommand\id{{\mathbf{i}}}
\newcommand\redn{{\mathrm{redn}}}

\newcommand\Mx{{\mathrm{M}}}
\newcommand\fr{{\mathbf{f}}}
\newcommand\D{{\mathbf{D}}}
\newcommand\E{{\mathbf{E}}}
\newcommand\A{{\mathbf{A}}}

\newcommand\ssll{{\mathfrak{sl}}}
\newtheorem{thm}{Theorem}[section]
\newtheorem{prop}[thm]{Proposition}
\newtheorem{obs}[thm]{Observation}
\newtheorem{lem}[thm]{Lemma}
\newtheorem{defn}{Definition}
\newtheorem{alt}{Equivalent Formulation}[section]
\begin{document}
\bibliographystyle{plain}
\title{Extrinsic properties of automorphism groups of formal groups}
\author{Jonathan D. Lubin and Ghassan Y. Sarkis}
\date{}
\maketitle
%
%
%
%\input extrinAbstract
%%%%%%% extrinAbstract.tex, for
%%%%%%%
%%%%%%%
%%%%%%%
\begin{abstract}
We prove two conjectures on the automorphism group of
a one-dimen\-sional formal group law defined over a field
of positive characteristic.
The first is that if a series commutes with a nontorsion
automorphism of the formal group law,
then that series is already an
automorphism. The second is that
%in the case of height greater than unity,
the group of automorphisms is its own normalizer in
the group of all invertible series over the ground field.
A consequence of these results is that a formal group
law in positive characteristic is determined by any one
of its nontorsion automorphisms.
\end{abstract}
%
%\input extrin010
%
%%%%%%% extrin010.tex, for the introduction of extrin.tex
%%%%%%%
%%%%%%%
%%%%%%%
\section*{Introduction}
In this paper we deal with formal group laws over a field
of characteristic $p>0$, in other words with coordinatized
one-dimensional formal groups. The theory of one-dimensional
formal groups in positive characteristic is rather different in
flavor from that in characteristic zero. In particular the
fact that a formal group
in characteristic zero is determined by one of its endomorphisms
is so immediate that it hardly deserves notice. From $f\in\End({G})$,
if ${G}$ is a formal group with coefficients lying in any $\QQ$-algebra,
one can derive by a simple degree-by-degree computation the
logarithm $\log_{G}\colon{G}\to\scrA$, where $\scrA$ is
the additive group law $x+y$, and from this the coefficients
of ${G}(x,y)$ drop out, again degree by degree.
\par
In characteristic $p$, however, the picture is entirely different:
for instance, for $x^p$ to be an endomorphism, it is only necessary
that ${G}(x,y)\in\FF_p[[x,y]]$. So if there is any hope of
characterizing a formal group law by one or more of its
endomorphisms, these should be automorphisms.
Further, a series such as $\alpha x\in k[[x]]$ will
be an automorphism of ${G}(x,y)$ if and only if
the only nonzero coefficients of ${G}$ occur in degrees
congruent to $1$ modulo $m-1$, for $m$ the multiplicative
period of $\alpha$. So we mostly worry about nontorsion
automorphisms of ${G}$.
\par
The first result along these lines was Theorem 7 of
\cite{gs:thesis}, and it states that if $u$ is a
nontorsion element of the center $\ZZ_p^*$ of the absolute
automorphism group of a formal group law ${G}$
and $\psi(x)$ is a series commuting with~$u$, then
$\psi$ is an endomorphism of ${G}$. A more recent one
is due to Li, who generalized the methods
of \cite{gs:thesis} in \cite{li:preprint}
to the case where $u'(0)$ generates $\FF_{p^h}$.
\par
The current paper's first main result is the strongest possible
theorem of this type: in the group $k[[x]]^\circ$ of
all invertible series over $k$, the centralizer of a
nontorsion automorphism of ${G}$ is contained in $\Aut({G})$.
This is certainly of interest in itself, but is also an
important ingredient in the proof of the other main result.
This theorem says that $\Aut({G})$ is its own normalizer
in~$k[[x]]^\circ$.
\par
And in turn, this theorem shows that, as a subgroup
of $k[[x]]^\circ$, a given automorphism
group comes from only one formal group law; this can then be
used to see that a single nontorsion automorphism
determines completely the formal group law that it
belongs to. One must hasten to add that there seems to be no
calculation of any effective nature that would allow computation
of any of the coefficients of this mysterious group law.
\par
%
%
%\input extrin020
%%%%%%% extrin020.tex, for the Ònotations and conventionsÓ of extrin.tex
%%%%%%%
%%%%%%%
%%%%%%%
\section{Notations, conventions, background}\label{sec1}
%
%\input extrin021 %(extrin021.tex)
%%%%%%% extrin021.tex, for the Ònotations and conventionsÓ of extrin.tex
%%%%%%%
%%%%%%%
%%%%%%%
Our formal group laws ${G}$ will all be of dimension one, defined over a
field $k\/$ of characteristic $p>0$, and of finite height greater
than $1$.
Since our interest throughout this paper is in coordinatized
objects, we will feel free to omit the word ``law'' and speak
somewhat inaccurately of ``formal groups''.
Indeed, the significance of the word ``extrinsic'' in our
title is that we are dealing with the group of automorphisms
of ${G}$ \textit{as series}, and how these groups may sit inside the
larger group of all invertible power series over~$k$.
\par
The one-dimensionality of our formal groups implies that they
are commutative, so that the set of endomorphisms of ${G}$, written
$\End_k({G})$, is a ring, the addition being ${G}$-addition of
series, and multiplication being composition.
The natural map $\ZZ\to\End_k({G})$ is denoted $n\mapsto[n]_{G}$;
and this map extends canonically to $\ZZ_p$, so that the endomorphism
ring has a natural structure of $\ZZ_p$-algebra.
Over any algebraically closed field $\Omega$ containing $k$,
it has been known since Dieudonn\'e \cite{dieu} that the endomorphism
ring is isomorphic to the maximal order in a central division
algebra over $\QQ_p$ of rank $h^2$ and invariant $1/h$, where
$h\/$ is the height of ${G}$.
When dealing with such a formal group ${G}$,
we will further assume of $k\/$ that this field is so large
that $\End_k({G})=\End_\Omega({G})$. Such a field $k$
must necessarily contain $\FF_{p^h}$, but this condition may
not be sufficient to catch all the endomorphisms of ${G}$.
\par
One may use ${G}(x,y)$ to add not only endomorphisms, but
arbitrary series in $xk[[x]]$: if $\varphi$ and $\psi$ are
two such, then their ${G}$-sum
is $(\varphi+_{G}\psi)(x)={G}(\varphi(x),\psi(x))$.
Using the additive $\ZZ$-valued valuation $v_x$ on the
ambient ring $k[[x]]$ and the shape of ${G}(x,y)$ as
$x+y+$ (higher terms), we see that if $v_x(\psi)>v_x(\varphi)$,
then $v_x(\varphi+_{G}\psi)=v_x(\varphi)$; for
similar reasons, if $v_x(\varphi)=v_x(\psi)$, then
$v_x(\varphi-\psi)=v_x(\varphi-_{G}\psi)$. Note that
the corresponding statement about the two kinds of addition
does not hold.
\par
We will need to make reference to the closeness of a series $u(x)$
to the identity series $x$, and in our context it is most
appropriate to use the nonstandard
measure of proximity $w(u)=v_x(u(x)-x)$. This instead of
$v_x\bigl(\frac{u(x)}{x}-1\bigr)$, which is usually more natural.
\par
As we have said above, the ring of $k\/$-endomorphisms of ${G}$
is denoted $\End_k({G})$; the group of
its invertible elements is denoted
$\Aut_k({G})$. Because of our assumption that $k\/$ is
sufficiently large, we will drop the subscript. Much of the
time, we will be wanting to consider the ring and its units
as abstract topological ring and group, so we denote them $\E$
and $\A$, respectively; then $\A=\E^*$. We will also be dealing
with the (noncommutative) fraction field of~$\E$,
$\D=\E\otimes_{\ZZ_p}\!\!\QQ_p$. We will have occasion
later to use the fact that there are algebraic extensions
$K$ of $\QQ_p$ such that $\D\otimes_{\QQ_p}\!K\cong\Mx_h(K)$,
the ring of $h$-by-$h$ matrices over $K$.
%the ring of $h$-by-$h$ matrices over $K$;
%that $\D\otimes_{\ZZ_p}\!\frO\cong\Mx_h(\frO)$, where
%$\frO$ is the ring of integers of $K$; and that extension
%of the base from $\ZZ_p$ to $\frO$ also turns $\A$ into something
%isomorphic to $\GL_h(\frO)$.
\par
If $f(x)\in k[[x]]$, we will denote its $n$-fold iterate $f^{\circ n}$,
but when $f\thinspace$ is in $\End({G})$ and is considered as element of
$\E$, then this iterate will be denoted simply $f^n$.
\par
We refer the reader in search of a readable introduction to
formal groups to the classic text \cite{fro}.
\par
%\input extrin030
%%%%%%%% extrin030.tex, for the second content section of extrin.tex
%%%%%%%
%%%%%%%
%%%%%%%
\section{Information from a single automorphism}\label{IFSA}
%
%\input extrin031 %(extrin031.tex)
%%%%%%%% extrin031.tex, in the first content section of extrin.tex
%%%%%%%
%%%%%%% called from extrin030.tex
%%%%%%%
How much can you tell about a formal group ${G}$ in characteristic
$p>0$ by looking at one of its automorphisms $u(x)$?
In the case of nontorsion automorphisms, there is one
piece of information that comes fairly easily and cheaply:
by looking at how fast certain iterates of $u$ approach
the identity, we can determine the height of ${G}$. Let us
explain this by going into the arithmetic of $\E$ a little more
deeply.
\par
The division ring $\D$ has on it an additive valuation extending
the canonical $\ZZ$-valued $v_p$ for which $v_p(p)=1$. Namely, for
$z\in\D$, $V(z)$ is defined
to be $v_p(\redn(z))/h$, where $\redn\colon\D\to\QQ_p$
is the reduced norm. Recall that this multiplicative homomorphism
is actually defined
%
%
% IÕm changing the following from a statement about $\QQ_p$-algebras
% to a general one about $K$-algebras for a field $K$.
%
%
for any central simple finite-dimensional $K$ algebra $A$,
and has the property that if $\dim_KA=h^2$
and $F$ is a commutative field with
$K\subset F\subset A$ and of degree $h$ over $K$,
then the reduced norm agrees on $F$ with the ordinary $F$-over-$K$
norm of Galois theory.
This is not a definition of the reduced norm, but
it is a sufficient description for our purposes. In case
$A$ is the ring of matrices over $K$, the reduced norm
is the determinant. One shows easily
that $V\colon\D\to\frac{1}{h}\ZZ\cup\{\infty\}$
is an additive valuation in
the usual sense: $V(zz')=V(z)+V(z')$ and $V(z+z')\ge\min(V(z),V(z'))$.
\par
Although it represents a break in the presentation, it is most
convenient to insert here a mention of the commutator of two
elements $a$ and $b$ of~$\E$, even of~$\D$.
We define $[a,b]$ to be $ab-ba$.
Then as may be easily verified, the commutator is
$\ZZ_p$-bilinear, $V([a,b])\ge V(a)+V(b)$, and if $V(a)$ and $V(b)$
both are positive, then $V([1+a,1+b])\ge V(a)+V(b)$. A little less
obvious is the following:
%
%%%%% PROPOSITION ON COMMUTATOR OF AN ARB. ELEMENT WITH A POWER
%%%%% OF SOMETHING IN THE MAXIMAL IDEAL
%
\begin{prop}\label{cmmttrprops}\textsl{
If $a, b\in\D$, then $V([a^{n+1},b])\ge V([a,b])+nV(a)$. Moreover, if
$V(a)>\frac{1}{p-1}\thinspace$,
then $V\Bigl([(1+a)^p,b]\Bigr)=1+V([a,b])=1+V([1+a,b])$.
}\end{prop}
\par
The first statement follows from the fact that
\begin{equation*}
[a^{n+1},b] = \sum_{i=0}^na^i[a,b]a^{n-i}\,.
\end{equation*}
The second follows from the fact that in the expansion
of $[(1+a)^p,b]$ coming from the binomial expansion,
the first term is zero, while the $V$-value of the second term
is definitely less than that of any of the others.
\par
The hypothesis that $h$ is the height of ${G}$ is exactly that
$v_x([p]_{G})=p^h$. Since $v_x(\varphi\circ\psi)=v_x(\varphi)v_x(\psi)$,
we conclude that for $z\in\End({G})$, $V(z)=\frac{1}{h}\log_p(v_x(z))$.
So we have:
\par
\begin{prop}\label{s3p1}\textsl{
Let $u$ be an auto\-morphism of ${G}$,
%%%%%%with $u(x)\equiv x\pmod{x^2}$,
and call the identity auto\-morphism $\id(x)=x$. Then
%$V(u-\id)=\frac{1}{h}\log_p(w(u))$.
\[
w(u)=v_x(u-_{G}\id)=p^{hV(u-1)}\,,
\]
where in the rightmost member of the display, $u$ is thought of
as an element of~$\E$.
}\end{prop}
\par
Any element of $\E$ sits in a commutative subfield, of
residue field extension degree dividing $h$, so that
a unit $u^{p^h-1}$ will always be congruent to $1$ modulo the
maximal ideal of $\E$. But once a quantity is this
close to the identity, we understand fully how rapidly its successive
$p$-th powers approach $1$, just from the expansion of $(1+\varepsilon)^p$:
\par
%
%\input extrin032 %(extrin032.tex)
%
%%%%%%% extrin032.tex, in the first content section of extrin.tex
%%%%%%%
%%%%%%% called from extrin030.tex
%%%%%%%
\begin{obs}\label{s3p2}\textsl{
If $F$ is an algebraic extension of $\QQ_p$ and $z\in F$ such that
$v_p(z-1)>0$, then:
\[
v_p(z^p-1)\begin{cases}
=pv(z-1)&\text{if $\thinspace v(z-1)<\frac{1}{p-1}$}\\
\ge p/(p-1)&\text{if $\thinspace v(z-1)=\frac{1}{p-1}$}\\
=1+v(z-1)&\text{if $\thinspace v(z-1)>\frac{1}{p-1}\,.$}
\end{cases}
\]
}\end{obs}
\par
By combining \ref{s3p1} and \ref{s3p2}, we conclude:
\begin{thm}\label{newtpol}\textsl{
Let $u$ be an auto\-morphism of the height-$h$ formal group ${G}$,
with $u(x)\equiv x\pmod{x^2}$. Then $w(u^{\circ p})$ depends on
$w(u)$ in the following way:
\[
w(u^{\circ p})
\begin{cases}
=w(u)^p&\text{if \thinspace $w(u)<p^{h/(p-1)}$}\\
\ge p\thinspace^{ph/(p-1)}&\text{if $\thinspace w(u)=p^{h/(p-1)}$}\\
=p^hw(u)&\text{if $\thinspace w(u)>p^{h/(p-1)}\,.$}
\end{cases}
\]
}
\end{thm}
\par
The doubt inherent in the second rule is unavoidable: the
only way to get around it would be to add
to the hypotheses an additional
statement about the proximity of $u$ to a $p$-th root of unity.
In any event, the directly observable quantity $w(u)$
and the perfectly computable corresponding quantities for the
$p$-power iterates of $u$ exhibit the height of ${G}$, since
$p^h$ is the stable value of $w(u^{\circ p^{m+1}})/w(u^{\circ p^m})$
as $m\to\infty$.
\par
\begin{defn}\textnormal{
When an automorphism $u$ of ${G}$ satisfies the
condition that $w(u)>p^{h/(p-1)}$, we say that
$u$ is \textsl{in the stable range}.
}\end{defn}
\par
%
%
%\input extrin040
%
%%%%%%% extrin040.tex, for the third content section of extrin.tex
%%%%%%%
%%%%%%%
%%%%%%%
\section{Series that commute with an automorphism}
%
%\input extrin041 %(extrin041.tex)
%
%%%%%%% extrin041.tex, in the first content section of extrin.tex
%%%%%%%
%%%%%%% (called from extrin040.tex )
%%%%%%%
The result of this section answers almost completely the question
of how a general $k$-series can commute with a ${G}$-automorphism.
The only question remaining open is whether a nonendomorphic series
can commute with a $p$\nobreakdash-power torsion automorphism, and we do not have
even partial results there.
\par
\begin{thm}\label{th1}\textsl{
Let ${G}$ be a formal group of finite height, defined over a
field $k$ of characteristic $p>0$, with a nontorsion automorphism
$u$. If $\psi(x)\in k[[x]]$ with $u\circ\psi=\psi\circ u$, then
$\psi\in\End({G})$.
}\end{thm}
\par
This theorem is new and of interest only in the case that $h>1$,
since in the height-one case, all endomorphisms are in $\ZZ_p$,
and Theorem 7 of \cite{gs:thesis} shows that a series commuting
with a nontorsion element of $\ZZ_p^*\subset\Aut({G})$ is an
endomorphism.
The proof requires a few new concepts and lemmas, which we now present.
The first remark to make is that since any two formal groups
of the same height are isomorphic over any algebraically closed
field containing $k$, it suffices to prove \ref{th1} for a single
formal group ${G}$ of each height. The most important standardization
that we
choose is of the kind mentioned in Section 3.3.2 of \cite{gs:thesis},
where the formal group ${G}\in\FF_p[[x,y]]$ satisfies
\[
{G}(x,y)=x+y+{G}_0\Bigl(x^{p^{h-1}},\,y^{p^{h-1}}\Bigr)
\quad\hbox{and}\quad[p]_{G}(x)=x^{p^h}+\cdots\,,
\]
but for our purposes, all we will be using
of the first condition is the fact that
the two partial derivatives of ${G}$ are constant $1$.
To construct a height-$h$ formal group
satisfying both conditions, one may start with, say, the
$p$-typical logarithm $x+\sum_{m>0}p^{-m}x^{p^{mh}}$ to get
a formal group law in characteristic zero, reduce it to
get a group law $F$ with $[p]_F(x)=x^{p^h}$, and then apply
Lemma 23 of \cite{gs:thesis} to get an $\FF_p$-isomorphic
formal group ${G}$ of the form $x+y+G_0(x^p,y^p)$,
for which the $p$-endomorphism will still be $x^{p^h}$.
For details of the construction of a formal group from
a logarithm in characteristic zero, refer to \cite{haze}.
Our formal group has one other property that we will be needing,
and that is that for every $\alpha\in\FF_{p^h}$ and every
$r$, there is an endomorphism whose lowest monomial is $\alpha x^{p^r}$.
This is a general fact about formal groups defined
over $\FF_p$, not dependent on our construction of~${G}$.
From the fact that $[p]_{G}(x)=x^{p^h}$, it follows that
all endomorphisms of ${G}$ are defined over~$\FF_{p^h}$.
\par
%
%\input extrin042 %(extrin042.tex)
%
%%%%%%% extrin042.tex, in the first content section of extrin.tex
%%%%%%%
%%%%%%% (called from extrin040.tex )
%%%%%%%
A consequence of this standardization of ${G}$ is that if
an automorphism $u$ has $u'(0)=1$, then $u'(x)=1$, constant.
For, $u-_{G}\id$
is an endomorphism without first-degree term, so that it is a power
series of the form $\gamma(x^p)$, and hence $u={G}(x,\gamma(x^p))$,
whose derivative is $1$.
\par
For uniformity, we will denote series not known to be endomorphisms
by Greek letters, endomorphisms by Latin. In particular, we will call
the Frobenius $\fr(x)=x^p$, an endomorphism because of our assumption
that ${G}$ is defined over the prime field. For any $\rho\in k[[x]]$,
we have $\fr^{\circ r}\circ\rho=\rho^{(p^r)}\circ\fr^{\circ r}$,
where the exponent on $\rho$ is the method we use to denote that
each coefficient of $\rho$ has been raised to the $p^r$-power.
\par
\begin{defn}\textnormal{
For series $\varphi$ and $\psi$ in $xk[[x]]$, the
\textsl{${G}$-commutator} of the two, $[\varphi,\psi]$, is
$\varphi\circ\psi-_{G}\psi\circ\varphi$.
}\end{defn}
\par
Since the two compositions have the same initial degree,
it follows that
$v_x([\varphi,\psi]) = v_x(\varphi\circ\psi-\psi\circ\varphi)$.
\par
\begin{lem}\label{LA}\textsl{
Let $u$ and $\psi$ be commuting $k$-series, with $u\in\End({G})$.
Then $\psi'(0)\in\FF_{p^h}$.
}\end{lem}
This proposition appears in \cite{gs:thesis}, but we repeat it here
for completeness' sake.
Starting with a general nontorsion $u$,
%this series may be replaced
%by $u^{\circ(p^h-1)}$, which still commutes with $\psi$ and has
%first-degree coefficient equal to $1$. Once more we replace
%this series by a suitable $p$-power iterate, using
%Theorem \ref{newtpol}, so that we now are working with
%$u$ satisfying $w(u)>p^{h/(p-1)}$.
we replace $u$ by a suitable iterate that is
in the stable range, i.e. $w(u)>p^{h/(p-1)}$, so that
$u$ falls into the last
case of Theorem \ref{newtpol}.
Throughout the rest of the proof of Theorem \ref{th1}, we
will assume that $u$ has this property.
\par
%
%\input extrin043 %(extrin043.tex)
%
%%%%%%% extrin043.tex, in the first content section of extrin.tex
%%%%%%%
%%%%%%% (called from extrin040.tex )
%%%%%%%
Now, for
$u(x)\equiv x + \lambda x^{p^r}\!\!\!\pmod{x^{p^r+1}}$, where $r>h/(p-1)$
and $\lambda\ne0$, we also have
$u^{\circ p}(x)\equiv x+\lambda'x^{p^{r+h}}\!\!\!\pmod{x^{p^{r+h}+1}}$.
Write $\psi(x)\equiv \alpha x\!\!\!\pmod{x^2}$, and
perform the two compositions of $u$ and $\psi$ to get the congruences
\begin{eqnarray*}
\psi(u(x))&\equiv&\psi(x)+\psi'(0)\lambda x^{p^r}
\\
u(\psi(x))&\equiv&\psi(x)+\lambda(\alpha x)^{p^r}\,,
\end{eqnarray*}
both modulo $(x^{p^r+1})$, which gives us the
equality $\lambda\alpha=\lambda\alpha^{p^r}$, so
that $\alpha\in\FF_{p^r}$. The same argument applied to $u^{\circ p}$ and
$\psi$ shows that $\alpha\in\FF_{p^{r+h}}$, and these two facts imply that
$\alpha\in\FF_{p^h}$.
\par
\smallskip
%\begin{defn}\textnormal{
%If $g$ is a ${G}$-endomorphism,
%$\C=\C_g\subset xk[[x]]$ is
%the set of all series
%$\rho$ such that $[g,\rho]\in\End({G})$.
%}\end{defn}
%\par
The ${G}$-commutator does not behave as well as one might like: in
particular, it is not bilinear, nor even biadditive, although one
checks easily that  if $g$ is an endomorphism, $[g,\rho]$ is
$\ZZ_p$-linear in $\rho$. As a result, under $G$-addition,
the set of all series $\rho$ such that $[g,\rho]\in\End({G})$
%$\C_g$
is a left-$\ZZ_p$-module containing both $\End({G})$ and all
series commuting with $g$.
\par
%\smallskip
Since $\End({G})$ is topologically closed in $xk[[x]]$, for any %%series
$\psi(x)\in xk[[x]]$ there is is an endomorphism $g$ whose distance from
$\psi$ is least---if $\psi\notin\End({G})$, there will be many. So,
assuming that $\psi$ is not an endomorphism, we will choose $g$ and $\delta$
such that $\psi=g+_{G}\delta$ and $v_x(\delta)$ is maximum.
From \ref{LA} it follows that $v_x(\delta)>1$,
since we may ${G}$-subtract from $\psi$
any endomorphism with the same first-degree monomial.
\par
Indeed, to prove Theorem \ref{th1} our aim will be to show that the
lowest monomial in $\delta$ is of the form $\alpha x^{p^r}$
with $\alpha\in\FF_{p^h}$, for if so, we would be able to subtract
an endomorphism with the same leading monomial from $\delta$ to get
a higher $v_x$\nobreakdash-value.
\par
%
%\input extrin044 %(extrin044.tex)
%%%%%%% extrin044.tex, in the first content section of extrin.tex
%%%%%%%
%%%%%%% (called from extrin040.tex )
%%%%%%%
First we show that the initial degree of $\delta$ is a power of $p$.
To do this, we first write $\delta(x)=\Delta(x^{p^r})$, with $\Delta'(x)\ne0$.
Now consider the ${G}$-endomorphism $z=[u,\delta]$, for which we have:
\begin{eqnarray*}
z&=&u\circ\delta-_{G}\,\delta\circ u\\
&=&u\circ\Delta\circ\fr^{\circ r}-_{G}\,\Delta\circ\fr^{\circ r}\circ u\\
&=&u\circ\Delta\circ\fr^{\circ r}-_{G}\,\Delta\circ u^{(p^r)}\circ\fr^{\circ r}\\
&=&\Bigl(u\circ\Delta-_{G}\,\Delta\circ u^{(p^r)}\Bigr)\circ\fr^{\circ r}\,,
\end{eqnarray*}
which at the very least shows that $z=z_0\circ\fr^{\circ r}$, for
an endomorphism $z_0=u\circ\Delta-_{G}\,\Delta\circ u^{(p^r)}$.
Recalling that $u(x)=x+\beta(x^p)$ and that $G(x,y)=x+y+G_0(x^p,y^p)$,
we take the
equality
\begin{equation}
z_0+_{G}\Delta\circ u^{(p^r)}=u\circ\Delta\tag{A}
\end{equation}
and conclude,
up to $p$-th powers, that $z_0+\Delta=\Delta$. That is, $z_0$ involves
only $p$-th powers, and so its derivative is zero. Differentiating
equation (A) and taking account of the shape of $u$ and $G$, we
see that $\Delta'\circ u^{(p^r)}=\Delta'$. But since $u^{(p^r)}$ is
a nontorsion invertible series, $\Delta'$ must be constant, nonzero
by construction, so that the initial degree of $\delta$ is indeed $p^r$,
and as we have observed, $r>0$.
\par
It will take a bit more work than above to show that the coefficient
of~$x^{p^r}$ in~$\delta$ is in $\FF_{p^h}$.
Let $u(x)=x+\lambda x^{p^m}+\cdots$ and, as
above, $\delta(x)=\alpha x^{p^r}+\cdots$. Our aim is to show
that the coefficient of $x^{p^{r+m}}$ in $[u,\delta]$
is the same as the coefficient
of $x^{p^{r+m+h}}$ in $[u^{\circ p},\delta]$. From an explicit
computation of what these coefficients are, we will be able
to conclude that $\alpha\in\FF_{p^h}$.
\par
From the definition of $\delta$ as fitting into the equation
$\psi=g+\delta$ and the hypothesis that $[u,\psi]=0$, we conclude
that $[u,\delta]=[g,u]$,
which is the commutator of two endomorphisms in
just the sense that we mentioned in Section \ref{IFSA}.
The automorphism $u$ is of the form $\id+_{G}a$ mentioned in
Proposition \ref{cmmttrprops}, so that we can say that
$v_x([u^{\circ p},\delta])=p^h v_x([u,\delta])$.
The direct computation that we are about to do shows that
$v_x([u,\delta])\ge p^{r+m}$. If greater, then the coefficients
we are interested in are both zero; if equal, then our special
normalization $[p]_{G}(x)=x^{p^h}$ shows again that
the coefficients are equal.
\par
%\input extrin045 %(extrin045.tex)
%%%%%%% extrin045.tex, in the first content section of extrin.tex
%%%%%%%
%%%%%%% (called from extrin040.tex )
%%%%%%%
Let us look first at the coefficient of $x^{p^{r+m}}$ in
$[u,\delta]$. Modulo terms of degree greater than $p^{r+m}$,
we have $u(\delta(x))\equiv\delta(x)+\lambda\alpha^{p^m}x^{p^{r+m}}$;
the other composition is
$\delta(u(x))=\Delta\Bigl((x+\lambda x^{p^m})^{p^r}\Bigr)
\equiv\Delta(x^{p^r})+\alpha\lambda^{p^r}x^{p^{r+h}}$, which
shows that $[u,\delta]
\equiv(\lambda\alpha^{p^m} - \alpha\lambda^{p^r})x^{p^{r+m}}
\pmod{x^{1+p^{r+m}}}$.
Now we use the hypothesis on $u$, that $w(u)>p^{h/(p-1)}$,
which guarantees that $u^{\circ p}$ has the form
$x+\lambda'x^{p^{r+h}}+\cdots$, but we will now use
the agreed-upon standardization of ${G}$ to show that $\lambda'=\lambda$.
Writing $u=\id+_{G}a$ for an endomorphism $a$ which
begins with the monomial $\lambda x^{p^r}$, and using the
binomial expansion and the fact that $v_x(a)>p^{h/(p-1)}$, we
see that $u^{\circ p}\equiv\id+_{G}a\circ[p]_{G}\pmod{x^{1+p^{r+h}}}$,
but our standardization is that $[p](x)=x^{p^h}$, so that
$u(x)=x+\lambda x^{p^{r+h}}+\cdots$, as claimed.
Now, precisely the same computation as before gives
$[u^{\circ p},\delta]
\equiv(\lambda\alpha^{p^{m+h}} - \alpha\lambda^{p^r})x^{p^{r+m+h}}
\pmod{x^{1+p^{r+m+h}}}$. The upshot is that we have
the equation
\begin{equation*}
\lambda\alpha^{p^{m+h}} - \alpha\lambda^{p^r}
=
\lambda\alpha^{p^m} - \alpha\lambda^{p^r}\,,
\end{equation*}
and as we have seen, this shows that $\alpha\in\FF_{p^h}$, and so
finishes the proof of Theorem~\ref{th1}.
\par
%
%
%\input extrin050
%%%%%%% extrin050.tex, for the last content section of extrin.tex
%%%%%%%
%%%%%%%
%%%%%%%
\section{Series that normalize the automorphism group}
%
%\input extrin051 %(extrin051.tex)
%%%%%%% extrin051.tex
%%%%%%%
%%%%%%% (called from extrin050.tex )
%%%%%%%
In this section we will be working with a series $\psi(x)\in k[[x]]$
such that for every $u\in\Aut({G})$, we
have $\psi\circ u\circ\psi^{-1}\in\Aut({G})$. It may be of interest
to note that conjugation by such a series
leaves $\End({G})$ invariant as well. Indeed, the center
of $\Aut(G)$ will be invariant, so we
let $g\in\End({G})$ and
$u=[a]_{G}$ for $a\in\ZZ_p^*$, making $u$ a nontorsion central
automorphism of ${G}$. Then $\psi\circ u\circ\psi^{-1}=[a']_{G}$,
a nontorsion central automorphism which commutes with
$\psi\circ g\circ\psi^{-1}$. Since this last commutes with a nontorsion
element of $\ZZ_p^*\subset\Aut({G})$, it is a ${G}$-endomorphism.
Our aim in this Section is to prove:
\par
\begin{thm}\label{bigTh}\textsl{
Let ${G}$ be a formal group of finite height
%$h>1$
defined over a field $k$ of characteristic $p>0$,
$k$ being large enough for
all endomorphisms of ${G}$ to be defined over $k$.
Then $\Aut_k({G})$ is its own normalizer in $k[[x]]^\circ$.
}\end{thm}
\par
An equivalent statement is:
\par
\begin{alt}\textsl{
Under the same hypotheses as Theorem \ref{bigTh},
if $\thinspace{F}$ is a formal group over $k$
with $\Aut_k({F})=\Aut_k({G})$, then
the two series ${G}(x,y)$ and ${F}(x,y)$
are the same.
}\end{alt}
\par
Every $\psi$ in the normalizer of $\Aut_k({G})$ corresponds to a formal
group ${G}^\psi = \psi({G}(\psi^{-1}x, \psi^{-1}y))$
such that $\Aut_k({G})=\Aut_k({G}^\psi)$.
Conversely, suppose ${F}$ is another formal group over $k$
such that $\Aut_k({G})=\Aut_k({F})$.
Theorem \ref{newtpol} insures ${G}$
and ${F}$ have the same height, and so they are isomorphic
over the algebraic closure of $k$; any such isomorphism
would normalize $\Aut_k({G})$.
\par
%
%
%%%%%%%% insertion by GYS:
%
%
If the height of the formal group $G$ is $1$, then
Theorem 4.1 can be painlessly proven via \cite{LMS} as follows.
In the notation of pages 59---60 of \cite{LMS},
$$
e([1+p]_G)
=\lim_{n\to\infty}\left(\frac{(p-1)v_x\left(\frac{[(1+p)^{p^n}]_G}{x}-1\right)}{p^{n+1}}\right)
=\left\{\begin{array}{rl}p-1&\mathrm{if\;}p>2\\2&\mathrm{if\:}p=2\end{array}\right.
$$
Let $A=\{[(1+p)^z]_G\colon z\in\ZZ_p\}$.
By Th\'eor\`eme 5.9 of \cite{LMS}, we know that the
separable normalizer of $A$
is an extension of a finite group of order dividing $e([1+p]_G)$,
by the group $A$.
But $\Aut({G})$ is such a group,
obviously contained in this separable normalizer.
We will therefore concern ourselves only with the case $h>1$.
\par
%
%
%%%%%%%% end og GYS insertion
%
%
\medskip
%
%\input extrin052 %(extrin052.tex)
%%%%%%% extrin052.tex
%%%%%%%
%%%%%%% (called from extrin050.tex )
%%%%%%%
For most of the rest of this Section, we will be considering the
ring $\End({G})$ abstractly, and so will denote it, its
group of units, and its fraction field by the letters $\E$, $\A$,
and $\D$, as in Section \ref{sec1}. The operation $\theta$ of
conjugation by the series $\psi$ that normalizes $\A$ in the
ambient group $k[[x]]^\circ$ induces an isometric isomorphism
of the group $\A$, and consequently
%(because $\ZZ_p$
%and $\QQ_p$ have no nontrivial automorphisms)
the corresponding
action $\bar\theta$ on the Lie algebra of $\A$, namely on $\D$ with the Lie
bracket $zw-wz$, is a $\QQ_p$-Lie-algebra automorphism.
\par
More than that: the kernel of the reduced norm in $\D$ is
contained in $\A$, and thus is a subgroup of $\A$, which we call $\A_0$.
We call $\A'$ the commutator subgroup of $\A$, so that
we have $\A\supset\A_0\supset\A'$, and these two algebraic subroups of
$\A$, being both of dimension $h^2-1$, have the same Lie algebra,
namely the kernel of the trace from $\D$ to $\QQ_p$, which we will
denote $\D_0$. The operation $\theta$ leaves $\A'$ invariant, so
induces an automorphism $\bar\theta_0$ of the $\QQ_p$-Lie algebra $\D_0$.
Upon extension of the base from $\QQ_p$ to an algebraically
closed field $\Omega$, $\D_0$ becomes isomorphic
to $\ssll_h(\Omega)$ and $\D$ becomes isomorphic to $\Mx_h(\Omega)$.
Now we know (for instance from Theorem 5, p. 283 of \cite{jake})
that any automorphism of the $\Omega$-Lie algebra $\ssll_h(\Omega)$
is of the form $z\mapsto AzA^{-1}$ for some $A\in\GL_h(\Omega)$ or
(in case $h>2$)
of the form $z\mapsto-Az^{{\mathrm{t}}}A^{-1}$, where $z^{{\mathrm{t}}}$
denotes the transpose of $z$ considered as a matrix.
The identity component of the automorphism group is thus of index
two in the full group of automorphisms, and
in particular,
the square of any automorphism of $\ssll_h(\Omega)$ is in this identity
component.
\par
Now for an element $y$ of $\D$, the specification that for
every $z\in\D_0$, $\bar\theta_0(z)y=yz$ amounts to a set of
$\QQ_p$-linear conditions on $y$. If these linear conditions
have a nonzero solution after extension of the base to $\Omega$,
they already do so over $\QQ_p$.
We therefore are ready to show:
\par
\begin{prop}\label{almost}\textsl{
Let ${G}$ be a formal group of finite height $h>1$ over a field
$k$ such that all endomorphisms of ${G}$ are defined over $k$.
If $\psi(x)\in k[[x]]^\circ$ and $\psi$ normalizes $\Aut({G})$,
then in case $h=2$, $\psi\in\Aut({G})$;
otherwise, $\psi^{\circ2}\in\Aut({G})$.
}\end{prop}
\par
%
%\input extrin053 %(extrin053.tex)
%%%%%%% extrin053.tex
%%%%%%%
%%%%%%% (called from extrin050.tex )
%%%%%%%
Let $\varphi$ be $\psi$ or $\psi^{\circ2}$ according as $h$ is
equal to or greater than $2$, and represent by the
letter $\theta$ the action of conjugation by
$\varphi$ on $\A$,
so that by the remarks preceding
the statement, there is a nonzero $y\in\D$ for
which $\bar\theta(z)=yzy^{-1}$ for all $z\in\D_0$, the
kernel of the trace from $\D$ to $\QQ_p$.
Since $\bar\theta(pz)=p\bar\theta(z)$ and $p$ commutes
with all elements of $\D$, we may take $y\in\E$.
\par
Now consider the mapping $z\mapsto y^{-1}\bar\theta(z)y$,
which is identity on the Lie algebra $\D_0$ so that
the automorphism of $\A$
given by $a\mapsto y^{-1}\theta(a)y$
is identity in a neighborhood of $1$ in $\A_0$.
Retranslating this fact to a statement about series,
we see that there is a ${G}$-endomorphism $g$
such that for any automorphism $u$ with trivial reduced norm
and for which $u$ is sufficiently close to the identity series,
$\varphi\circ u\circ\varphi^{-1}\circ g=g\circ u$. That is,
$\varphi^{-1}\circ g$ commutes with the nontorsion
automorphism $u$.
But Theorem \ref{th1} now applies, showing that $\varphi^{-1}\circ g$
is an endomorphism of ${G}$. Using the fact that $g$ is a nonzero
endomorphism, we see that $\varphi^{-1}$ is an endomorphism, and an
automorphism because invertible.
\par
\medskip
%
%\input extrin054 %(extrin054.tex)
%%%%%%% extrin054.tex
%%%%%%%
%%%%%%% (called from extrin050.tex )
%%%%%%%
The only task remaining is to show that if $h>2$, then
not only $\psi^{\circ2}$ but the
series $\psi$ of Proposition \ref{almost} is actually an
automorphism. We should say here that according to a
communication from G. Prasad \cite{prasad:email}, which
appeals to Proposition 3, P. 226 of his paper \cite{prasad:paper},
any automorphism of $\A_0$ is a Lie-group automorphism, and
so is inner. Appeal to Theorem \ref{th1} in the way we have
already done would prove our desired result directly.
%We prefer,
%however, to continue in the more primitive vein already established
%here.
In order to make our presentation as self-contained as possible,
however, we will continue in the vein already established.
\par
\begin{lem}\label{firstcoeff}\textsl{
Let $\psi$ be in the normalizer in $k[[x]]^\circ$ of $\Aut({G})$, where
${G}$ is a formal group of finite height $h>1$ defined over $k$.
Then $\psi'(0)\in\FF_{p^h}$.
}\end{lem}
\par
In case $\psi$ is already an automorphism, the conclusion is true; so
the only case of concern is the one where $\psi^{\circ2}$ is an
automorphism but $\psi$ is not.
We will refer to such series as \textsl{exceptional}.
Conjugation by $\psi$ will leave
stable the subgroup $1+p\ZZ_p$ of the center $\ZZ_p^*$
of $\Aut({G})$. But the automorphism group of $1+p\ZZ_p$ is
$\ZZ_p^*$, so that the only involutory automorphism of $1+p\ZZ_p$
is $u\mapsto u^{-1}$. Thus we have $\psi\circ u\circ\psi^{-1}=u^{-1}$
for every $u$ in the subgroup of $\Aut({G})$ corresponding
to $1+\ZZ_p$. Suppose now that $u$ is in the stable range
(automatic if $p>2$),
and that $u(x)\equiv x+\lambda x^{p^{mh}}\pmod{x^{1+p^{mh}}}$
for some $m$ and some nonzero $\lambda$.
Then $u^{-1}(x)\equiv x-\lambda x^{p^{mh}}\pmod{x^{1+p^{mh}}}$
and if we call
$\alpha=\psi'(0)$ we have, in a way that is by now very familiar:
\begin{eqnarray*}
\psi(u(x))&\equiv&\psi(x)+\alpha\lambda x^{p^{mh}}
\\
u^{-1}(\psi(x))&\equiv&\psi(x)-\lambda(\alpha x)^{p^{mh}}\,,
\end{eqnarray*}
both congruences modulo $(x^{1+p^{mh}})$. So $\alpha=-\alpha^{p^{mh}}$,
and by using the same argument with $u^{\circ p}$ instead of $u$,
we get $\alpha=-\alpha^{p^{(m+1)h}}$, again implying that
$\alpha\in\FF_{p^h}$.
\par
We can now prove Theorem \ref{bigTh} in the case that $p>2$. Let
$\psi\in k[[x]]^\circ$, normalizing $\Aut_k({G})$, and let
$\psi'(0)=\alpha\in\FF_{p^h}$, in accordance with \ref{firstcoeff}.
There is an automorphism $U$ with $U'(0)=\alpha$, so we set
$\rho=U^{-1}\circ\psi$. But for any series $\rho$ with first-degree
coefficient equal to $1$, $\rho^{\circ b}$ makes sense for any
$b\in\ZZ_p$, and in particular for $b=1/2$. And
where $\rho$ is in the normalizer of $\Aut({G})$,
all these ``iterates'' of $\rho$ are in that group as well.
It follows from Proposition \ref{almost} that $\rho\in\Aut({G})$,
and from this that $\psi\in\Aut({G})$.
\par
\medskip
The case of characteristic two is rather more difficult to
handle by these methods.
\par
%
%
%\input extrin060
%%%%%%% extrin060.tex, for the last content section of extrin.tex
%%%%%%%
%%%%%%%
%%%%%%%
\section{Nonexistence of exceptional series in characteristic two}
%
%\input extrin061 %% (extrin061.tex)
%%%%%%% extrin061.tex
%%%%%%%
%%%%%%% (called from extrin060.tex )
%%%%%%%
%
An exceptional series is a $\psi(x)\in k[[x]]^\circ$
that normalizes $\Aut({G})$ but is not itself an
automorphism. We know that in the normalizer, the
group $\Aut({G})$ is a subgroup of index at most two,
and that if $\psi$ is exceptional, then $\psi'(0)\in\FF_{p^h}$.
%
%We can easily see from what has gone before that an exceptional
%series $\psi$ must be torsion: since it commutes with the
%automorphism $\psi^{\circ2}$ yet is not an automorphism,
%it must be that $\psi^{\circ2}$ is torsion, by \ref{th1}.
%
Consider an exceptional series $\psi$: it commutes with the
automorphism $\psi^{\circ2}$ yet is itself not an automorphism,
a contradiction to Theorem \ref{th1} unless $\psi^{\circ2}$
is torsion. Thus $\psi$ itself must be a torsion series.
\par
We now specialize to the case that $p=2$. If $\psi$ is
exceptional and of period $2^nm$ with $m$ odd, then
by replacing $\psi$ by $\psi^{\circ m}$, we will get
an exceptional series of $2$-power period, and with
first degree coefficient equal to $1$. We will
call such a series \textsl{$2$-exceptional}.
Let us note first what the period of such a series
can be.
\par
\begin{lem}\textsl{
The period of a $2$-exceptional series is greater than $4$.
}\end{lem}
\par
Let us show first that an exceptional series $\psi$ can not
be an involution. If it were, then conjugation by $\psi$
would induce a nontrivial involution $\bar\theta$ on the Lie
algebra $\D_0$; if this has so much as a single eigenvalue
equal to $1$, we can get elements of $\A_0$ close to the
identity and commuting with $\psi$, so that by Theorem \ref{th1},
$\psi$ would be an automorphism, contrary to our assumption
of exceptionality. So all eigenvalues would be $-1$, which means
that the involution $\bar\theta$ is the negative of the identity.
But this is not a Lie-algebra homomorphism, since $\D_0$ is not
commutative.
\par
The proof that $\psi$ can not be of period $4$ is similar, since
$\psi^{\circ2}$ would be an involutory automorphism of ${G}$.
But there is only one such, namely $[-1]_{G}$, which induces
identity on the Lie algebra $\D_0$. This means again that
$\bar\theta$ is an involutory automorphism of $\D_0$; again its
only possible eigenvalues are $\pm1$, and the rest of the
argument is the same.
\par
\begin{lem}\textsl{
If $\psi(x)$ is a $2$-exceptional series for the formal
group ${G}$ of height $h$, then $h$ is even, and
$\psi(x)\equiv x+\lambda x^r\pmod{x^{r+1}}$, with
$\lambda\ne0$ and $r<2^{h/2}$.
}\end{lem}
\par
Since $\psi$ is torsion of period at least $8$, some
even iterate of $\psi$ is an automorphism $g$ of ${G}$
of period four. Such an automorphism generates a
commutative field extension of $\QQ_2$ of
degree $2$, so that $2|h$. But we know that $v_p(g-1)=1/2$,
so that $w(g)=2^{h/2}$, in the notation of Section 2. Since
a $2$-power iterate of $\psi$ is equal to $g$, we must have
$w(\psi)<w(g)$.
\par
%
%\input extrin062 %% (extrin062.tex)
%%%%%%% extrin062.tex
%%%%%%%
%%%%%%% (called from extrin060.tex )
%%%%%%%
%
\begin{lem}\textsl{
Let $\psi$ be a $2$-exceptional series for the formal group ${G}$
of height $h>1$. Then $\psi(x)=x+_{G}\Delta(x^{2^m})$ for
a series $\Delta$ with $\Delta'(x)=\Delta'(0)\ne0$,
and $m\ne0$.
}\end{lem}
\par
We assume at the outset, as we may, that ${G}$ is parametrized
so that ${G}(x,y)=x+y+{G}_0(x^{2^{h-1}},y^{2^{h-1}})$,
with ${G}_0(x,y)\in\FF_2[[x,y]]$.
Now write $\delta=\psi-_G\id$ and $\delta(x)=\Delta(x^{2^m})$
with $\Delta'(x)\ne0$ and $m>0$.
From the preceding Lemma, we know that
$m < h/2$.
Let $u$ be a nontrivial ${G}$-automorphism in $1+4\ZZ_2$, and set
$v_p(u-1)=r$. Calling $z=u-u^{-1}$,
we then have $v_p(z)=r+1$. Now use the
equation $\psi\circ u=u^{-1}\circ\psi$, substituting
$\psi=\id+_{G}\delta$. So we have:
\begin{eqnarray*}
u(x)+_{G}\Delta\bigl(u(x)^{2^m}\bigr)
&=&
u^{-1}(x)+_{G}u^{-1}\bigl(\Delta\bigl(x^{2^m}\bigr)\bigr)
\\
z(x)\;=\;g(x^{2^{h(r+1)}})
&=&
u^{-1}\bigl(\Delta\bigl(x^{2^m}\bigr)\bigr)
-_{G}
\Delta\bigl(u\bigl(x^{2^m}\bigl)\bigr)\,,
\
\end{eqnarray*}
using the fact that $u\in\FF_2[[x]]$.
Since $m<h(r+1)$, we may write
$g(x^{2^t})=u^{-1}\circ\Delta-_{G}\Delta\circ u$
with $t$ some positive number.
In particular, the left-hand side of this
equation has derivative zero, while the
right has derivative $\Delta'-\Delta'\circ u$,
because the special parametrization of ${G}$
guarantees that $u'(x)=1$. It follows that
$\Delta'$ is a constant, nonzero by construction.
\par
From this Lemma, we see that the first term in
$\psi(x)$ after the linear one will be of $2$-power
degree. What remains is to show that the coefficient
of this next monomial is in $\FF_{2^h}$.
\par
\begin{lem}\textsl{
Let $\psi(x)\equiv x+c x^{2^m}\pmod{x^{1+2^m}}$,
a $2$-exceptional series for the formal group ${G}$
of height $h>1$. Then $c\in\FF_{2^h}$.
}\end{lem}
\par
The preceding Lemmas have shown that a $2$-exceptional
series does necessarily have the specified form, and
that $m<h/2$. We now take any nontrivial $u$ in
$1+4\ZZ_2\subset\Aut({G})$, and
write $v_p(u-1)=r\ge2$, so that the first terms of $u(x)$
are $x+x^{2^{rh}}$. This $u$ is in the stable range.
From the equation $\psi\circ u=u^{-1}\circ\psi$ we
get $\psi\circ u=u^{\circ(-2)}\circ u\circ\psi$;
in other words, using the
congruence $u^{\circ2}\equiv x\pmod{x^{2^{(r+1)h}}}$,
we can conclude that $\psi$ and $u$ commute modulo
degree $2^{(r+1)h}$.
\par
Using the notation of the preceding Lemma, where
$\psi=\id+_{G}\delta$ and $\delta(x)=\Delta(x^{2^m})$,
we conclude that $u$ and $\delta$ commute modulo
degree $2^{(r+1)h}$, while $u$ and $\Delta$
commute modulo degree $2^{(r+1)h-m}$. Now we
have the two congruences
$u(\Delta(x))\equiv\Delta(x)+c^{2^{rh}}x^{2^{rh}}$
and
$\Delta(u(x))\equiv\Delta(x)+\Delta'(0)x^{2^{rh}}$,
both modulo degree $2^{(r+1)h-m}$. Since $\Delta'(0)=c$
and $m<h/2$, we derive the relation $c=c^{2^{rh}}$.
By using the same argument with $u^{\circ2}$ instead of $u$,
we get $c=c^{2^{(r+1)h}}$, which tells us that
$c\in\FF_{2^h}$. This concludes the
proof of the Lemma.
\par
It now requires only a few words to complete the
proof of Theorem \ref{bigTh} in case the characteristic
is $2$. If there were any exceptional series, there would
be one that was closest to the identity, necessarily
$2$-exceptional, of the form
$\psi(x)=x+c x^{2^m}+\cdots$. But there is
also an automorphism $U$ that has the same first
two monomials, and $U^{-1}\circ\psi$ would be an
exceptional series closer to the identity than $\psi$,
which gives the necessary contradiction.
\par
\bibliography{extrin}

\begin{thebibliography}{1}

\bibitem{dieu}
Jean Dieudonn\'e.
\newblock Groupes de lie et hyperalg\`ebres de lie sur un corps de
  caract\'eristique $p>0$: V{I}{I}.
\newblock {\em Math. Annalen}, 134:114--133, 1957.

\bibitem{fro}
A.~Fr\"ohlich.
\newblock {\em Formal Groups}, volume~74 of {\em Lecture Notes in Mathematics}.
\newblock Springer, 1968.

\bibitem{haze}
Michiel Hazewinkel.
\newblock {\em Formal Groups and Applications}.
\newblock Academic Press, 1978.

\bibitem{jake}
Nathan Jacobson.
\newblock {\em Lie Algebras}.
\newblock Dover Publications, 1979.

\bibitem{LMS}
Fran{\c{c}}ois Laubie, Abbas Movahhedi, and Alain Salinier.
\newblock Syst{\`e}mes dymaniques non archim{\'e}diens et corps des normes.
\newblock {\em Compositio Math.}, 132:57--98, 2002.

\bibitem{li:preprint}
Hua-Chieh Li.
\newblock Lubin's conjecture on $p$-adic dynamical systems.
\newblock Preprint, 2004.

\bibitem{prasad:paper}
Gopal Prasad.
\newblock Triviality of certain automorphisms of semi-simple groups over local
  fields.
\newblock {\em Math. Annalen}, 218:219--227, 1975.

\bibitem{prasad:email}
Gopal Prasad.
\newblock untitled.
\newblock e-mail to Dinakar Ramakrishnan, July 2005.

\bibitem{gs:thesis}
Ghassan~Y. Sarkis.
\newblock On lifting commutative dynamical systems.
\newblock {\em Journal of Algebra}, 293:130--154, 2005.

\end{thebibliography}
\end{document}